\def\magnification{\afterassignment\m@g\count@}
\def\m@g{\mag\count@
  \hsize6.5truein\vsize8.9truein\dimen\footins8truein}

\mag=\magstephalf  
\documentclass{amsart}
\newtheorem{theorem}{Theorem}[section]
\newtheorem{lemma}[theorem]{Lemma}
\newtheorem{corollary}[theorem]{Corollary}

\newtheorem{proposition}[theorem]{Proposition}

\newcommand{\R}{{\Bbb R}}
\newcommand{\N}{{\Bbb N}}

\newcounter{rmnum}

\newenvironment{alphanum}{\begin{list}{{\rm (\alph{rmnum})}}
{\usecounter{rmnum}\def\makelabel##1{\hss\llap{##1}}
\setlength{\leftmargin}{37pt}\setlength{\itemindent}{0pt}
\setlength{\topsep}{5pt}\setlength{\parsep}{0pt}\setlength
{\itemsep}{0pt}}}{\end{list}}

\begin{document}
\title{Continuous selections and finite $C$-spaces}
\author{Vesko Valov\\[3pt] }
\address{Department of Mathematics,
University of Swaziland, Pr. Bag 4,
Kwaluseni, Swaziland
\linebreak
Current address: Nipissing University, 100 College Drive, North Bay,
ON P1B 8L7, Canada}

\email{valov@realnet.co.sz}


\begin{abstract}
Characterizations of paracompact finite $C$-spaces via continuous
selections avoiding $Z_{\sigma}$-sets are given. We apply these
results to obtain some properties of finite $C$-spaces.
Factorization theorems and a completion theorem
for finite $C$-spaces are also proved.
\end{abstract}

\maketitle

\section*{1. Introduction}
\setcounter{section}{1}
\setcounter{theorem}{0}

\bigskip
Finite $C$-spaces were recently introduced by Borst \cite{b1:99}
and classified with the help of transfinite dimension function
$dim_C$ (Borst considered only separable metric spaces, but his
definition can be extended for any normal space, see \cite{cha:99}).
Every finite $C$-space has property $C$ in the sense of
Addis and Gresham \cite{ag:78} and, in the realm of compact spaces,
both classes coincide. On the other hand, every finite $C$-space
is weakly infinite-dimensional in the sense of Smirnov (abbr., $S$-wid),
while $C$-spaces are weakly infinte-dimensional in the sense of
Alexandrov (abbr., $A$-wid). 
One of the main problems in infinite dimension theory is whether
every weakly infinite-dimensional compact metric space has property
$C$ (it is well known that both properties $S$-wid and $A$-wid are
equivalent for compact spaces). The positive solution of this problem
would imply
that finite $C$-spaces and $S$-wid spaces constitute the same class.
But perhaps the most
important implication from the positive solution of this problem
is a positive solution
of another old problem in dimension theory (according to \cite{cha:99}
and \cite{bp:99}, Pasynkov was who first asked this question in 1972-73):
Is any product of two weakly infinite-dimensional compact metric
spaces weakly infinite-dimensional?

In the present paper we give other two characterizations
of paracompact finite $C$-spaces in the terms of continuous
selections for set-valued maps. The first one is similar to
Uspenskij's description of paracompact $C$-spaces
\cite[Theorem 1.3]{vu:98} and to the \v{S}\v{c}epin and N. Brodskii theorem
\cite{sb:96}.
To state this result let us agree
some notations. All spaces are assumed to be Tychonoff.
A set-valued map $\theta :X\rightarrow 2^Y$,
where $2^Y$ is the family of nonempty subsets of $Y$, is called
strongly lower semi-continuous (br., strongly lsc) if for
every compact $K\subset Y$ the set
$\{x\in X: K\subset\theta (x)\}$ is open in $X$; $\theta$ is
said to be $m$-aspherical if all $\theta (x)$ are $C^m$-sets (recall
that a $A\subset B$ is $C^m$-embedded, or $m$-embedded in $B$,
if every continuous image of a
$k$-sphere in $A$, $k\leq m$, is contractible in $B$; when $A$ is
$m$-embedded in itself, then it is called $C^m$).
We also say that a sequence of maps $\theta _n:X\rightarrow 2^Y$
is increasing (resp., decreasing) if
$\theta _n(x)\subset\theta _{n+1}(x)$ (resp., 
$\theta _{n+1}(x)\subset\theta _n(x)$) for all $x\in X$ and
$n\in\N$; the sequence $\{\theta _n\}$ is aspherical provided
it is increasing and each 
$\theta _n(x)$, $x\in X$, is $C^{n-1}$-embedded in 
$\theta _{n+1}(x)$.

\begin{theorem}
For a paracompact space $X$ the following are equivalent:

\begin{alphanum}
\item $X$ is a finite $C$-space; 

\item For any space $Y$ and any aspherical 
sequence of strongly lsc maps
$\phi _n\colon X\to 2^Y$, there exists $m$ such that
$\phi _m$ admits a continuous selection;

\item For any space $Y$ and any aspherical sequence
$\phi _n\colon X\to 2^Y$ of set-valued open-graph maps
there exists
$m$ such that $\phi _m$ admits a continuous selection. 
\end{alphanum}
\end{theorem}
 
Another characterization of paracompact finite $C$-spaces is the
following selection theorem of Michael's type (see \cite[Theorem 1.1]{gv:99}
for a similar characterization of paracompact $C$-spaces). 

\begin{theorem}
For a paracompact space $X$ the following are equivalent:

\begin{alphanum}
\item  $X$ is a finite $C$-space;

\item Suppose that $Y$ is a Banach space and 
$\phi\colon X\to {\mathcal F}_c(Y)$ is an lsc map. Then, for every
decreasing $Z$-sequence $\{\psi _n\}$ for $\phi$ consisting of
closed-graph maps $\psi _n\colon X\to {\mathcal F}(Y)$,   
there exist $m$ and a
continuous selection for $\phi$ avoiding the map $\psi _m$;

\item Let $Y$ be a Banach space and $\phi\colon X\to {\mathcal F}_c(Y)$
be an lsc map. Then, for every decreasing $Z$-sequence 
for $\phi$ consisting of 
closed sets $F_n\subset Y$, there exists a continuous
selection for $\phi$ avoiding some $F_m$.
\end{alphanum}
\end{theorem}

\bigskip
Here, ${\mathcal F}(Y)$ stands for the closed sets $F\subset Y$ and
${\mathcal F}_c(Y)$ denotes the convex elements of ${\mathcal F}(Y)$
(the same
notation will be used when $Y$ is a subset of a vector space).  A set-valued
map $\theta\colon X\to 2^Y$ is lower semi-continuous (br., lsc) if
$\theta ^{-1}(U)=\{x\in X: \theta (x)\cap U\neq\emptyset\}$ is open in $X$ for
every open $U\subset Y$. The term "continuous selections avoiding a set" is
due to E. Michael \cite{m2:88}. If $\phi\colon X\to 2^Y$ is a set-valued map
and $F\subset Y$, then a map $f\colon X\to Y$ is a selection for
$\phi$ avoiding $F$ if $f(x)\in\phi (x)\backslash F$ for every $x\in X$.
Similarly, if $\psi\colon X\to 2^Y$ is another set-valued map, we say that
$f$ is a selection for $\phi$ avoiding $\psi$ provided
$f(x)\in\phi (x)\backslash\psi (x)$ for any $x\in X$.
We adopt the following definition for a $Z_n$-set: 
if $n\geq 0$ we say that
$F\in{\mathcal F}(Y)$ is a $Z_n$-set in $Y$ if
the set
$C(I^n,Y\backslash F)$
is dense in $C(I^n,Y)$,
where $C(I^n,Y)$ stands for the space of all continuous maps from $I^n$
into $Y$ equipped with the compact open topology.  
When $F=\cup\{F_n:n\in\N\}$ such that each $F_n$ is a $Z_n$-set
in $Y$, we say that it is a $Z_{\sigma}$-set in $Y$.
The collection
of all $Z_n$-sets in $Y$ is denoted by ${\mathcal Z}_n(Y)$. Obviously,
$Z_0$-sets are precisely the closed and nowhere dense subsets of $Y$.
A sequence of maps $\psi _n\colon X\to 2^Y$
is called a $Z$-sequence for a given map $\phi\colon X\to 2^Y$ if
$\phi (x)\cap\psi_n(x)\in{\mathcal Z}_{n-1}(\phi (x))$ for every $x\in X$
and $n\in\N$.

The definition and some general properties of $C$-spaces and
finite $C$-spaces are given in Section 2. 
Theorem 1.1 and Theorem 1.2 are proved in Section 3 and Section 4,
respectively. Section 5 is devoted to some applications. 
The idea to characterize dimension like 
properties by removing maps from "small" sets goes back
to the Alexandroff characterization of $n$-dimensional compacta
in Euclidean spaces $\R^k$ as the compacta which are removable 
from any $k-n-1$-dimensional polyhedron \cite{pa:28}, \cite{re:95}
(see \cite{drs:98} for more general treatment
of removable maps in Euclidean spaces). Following this general idea, we apply
Theorem 1.2 to characterize  
finite $C$-spaces in terms of removable sequences of maps (see Proposition 5.1). 
This characterization is then used to obtain that finite $C$-space
property is preserved by approximately invertible maps (an
analogue to Ancel's result \cite{fa:85}).
In the final
Section 6 we prove a factorization theorem for finite $C$-spaces
when the image space is metrizable and apply this theorem to
obtain that any metrizable finite $C$-space has a completion
with the same property. 

I cordially thank Prof. V. Gutev for his suggestion to redefine
$C$-space and finite $C$-space properties by using functionally open covers.
My thanks are also due to the referee for his/her valuable comments and to
Prof. T. Banakh for the preprint of his joint paper
\cite{bt:99} with K. Trushchak.

\section*{2. Some preliminary results}
\setcounter{section}{2}
\setcounter{theorem}{0}

\bigskip
This section is devoted to the definition and general properties of
$C$-spaces and finite $C$-spaces.

$C$-space property was introduced by Haver \cite{ha:74} for compact
metric spaces. Addis and Gresham \cite{ag:78} reformulated
Haver's definition so that it has meaning for any space:
A space $X$ has property $C$ (or $X$ is a $C$-space) if
for any sequence $\{\omega _n\}$ of open covers of $X$ there exists a sequence
$\{\gamma _n\}$ of open disjoint families in $X$ such that each $\gamma _n$
refines $\omega _n$ and $\bigcup\{\gamma _n:n\in\N\}$ covers $X$. 
The family $\{\gamma _n\}$ is called $C$-refinement for
$\{\omega _n\}$. We now adopt the following definition:
$X$ is a $C$-space if every sequence $\{\omega _n\}$ of
locally finite and functionally open covers of $X$ has a
$C$-refinement $\{\gamma _n\}$ such that 
$\bigcup\{\gamma _n:n\in\N\}$ is a locally finite and functionally 
open cover of $X$. In the realm of paracompact spaces, this
definition coincides with the Addis and Gresham one.
Let us note that countable-dimensional metric spaces have property
$C$ \cite{ag:78}. Since the
countable sum theorem for property $C$ holds in the class of paracompact
spaces \cite{gv:99}, any strongly countable-dimensional (a countable
union of closed finite-dimensional subspaces) paracompact is a $C$-space.
On the other hand, there exists a metric $C$-compactum which is not
countable-dimensional \cite{p:81}.

In \cite{b1:99} P. Borst originally defined finite $C$ property for
separable metric
spaces: $X$ is a finite $C$-space if for
any sequence $\{\omega _n\}$ of finite open covers of $X$
there exists a finite sequence $\{\gamma _n\}_{n=1}^{k}$ of disjoint
open families in $X$ such that each $\gamma _n$ refines
$\omega _n$ and $\bigcup\{\gamma _n:n=1,..,k\}$ covers $X$. 
The sequence $\{\gamma _n\}_{n=1}^{k}$ is called a finite 
$C$-refinement for $\{\omega _n\}$. Obviously, the existence of
a finite $C$-refinement for $\{\omega _n\}$
is equivalent to the existence of a finite $C$-refinement consisting of
finite families.
As above, we reformulate this definition for arbitrary spaces: $X$ is a finite
$C$-space if any sequence
$\{\omega _n\}$ of finite functionally open covers of $X$
has a finite $C$-refinement $\{\gamma _n\}_{n=1}^{k}$ 
such that each $\gamma _n$ is a finite functionally open 
and disjoint family.
Note that, in the realm of normal spaces
we can omit the requirement $\omega _n$ and $\gamma _n$ to be functionally
open.

Since any locally finite family of open sets in a pseudocompact
space is finite, we have the following

\begin{proposition}
Property $C$ and finite $C$-space property are equivalent for
any pseudocompact space.
\end{proposition}
 
The proof of our first observation follows from the well known fact that
if $Y\subset X$ is $C^*$-embedded (i.e., every bounded function on $Y$
can be continuously extended to a bounded function on $X$), then for any finite functionally open cover $\omega$ of $Y$ 
there exists a finite functionally open cover $\gamma$ of 
$X$ such that $\omega$ is refined by $\gamma\cap Y$.

\begin{proposition} For any space $X$ we have:
\begin{alphanum}
\item
If $X$ is a finite $C$-space and $Y\subset X$ is $C^*$-embedded,  
then $Y$ has finite $C$-space property;
\item $X$ is a finite $C$-space if and only if $\beta X$ is 
a $C$-space.
\end{alphanum}
\end{proposition}

Let us note that, in case $X$ is normal,
Proposition 2.2(b) was proved by Chatyrko \cite{cha:99} (he
proved somewhat more, that $dim_CX=dim_C\beta X$). 
Next lemma is well known in case of metrizable spaces (see \cite{b1:99}).
The same proof remains valid in our situation.

\begin{lemma}
Every normal finite $C$-space is $S$-weakly infinite dimensional.
\end{lemma}

Recall that a space $X$ satisfies the condition $(K)$ \cite{ep:86}
if there is a compact set $K\subset X$ such that every closed
in $X$ set which is disjoint from $K$ has finite covering
dimension $dim$. P. Borst \cite[Theorem 3.8]{b1:99} proved that
his definition of finite $C$-spaces is equivalent to the
following one: Any sequence of open covers has a finite 
$C$-refinement. This implies that, in the realm of metric spaces,
every finite $C$-space has property $C$. Next theorem shows
that the same is true for any space.  

\begin{theorem}
For a space $X$ the following condition are equivalent:
\begin{alphanum}
\item $X$ is a finite $C$-space;
\item Any sequence $\{\omega _n\}$ of locally finite
functionally open covers of $X$ has a finite $C$-refinement
$\{\gamma _n\}_{n=1}^{s}$ such that 
$\{\gamma _n:n=1,2,..,s\}$ is a locally finite 
functionally open cover of $X$.
\end{alphanum}
\end{theorem}

\begin{proof} $(a)\Rightarrow (b)$. Let
$\omega _n=\{U_{\alpha}^n:\alpha\in A_n\}$ be
a sequence of locally finite functionally open covers of $X$.
Then there exists a map $f\colon X\to Y$ onto a metrizable space
$Y$ and a sequence $\eta _n=\{W_{\alpha}^n:\alpha\in A_n\}$
of open covers of $Y$ such that
$f^{-1}(W_{\alpha}^n)\subset U_{\alpha}^n$ for every $n$ and
$\alpha\in A_n$ (see \cite[Exercise 5.1.J(b)]{re:89}). Let
$Z=(\beta f)^{-1}(Y)$, where $\beta f\colon\beta X\to\beta Y$ is
the extension of $f$, and $\overline{\omega _n}=(\beta f)^{-1}(\eta _n)$.
Then each $\overline{\omega _n}$ is a functionally open cover of $Z$
and its restriction on $X$ refines $\omega _n$. By Proposition 2.2(b),
$Z$ is a finite $C$-space. Hence, by Lemma 2.3, $Z$ is $S$-wid.
Because every paracompact $S$-wid space satisfies the condition $(K)$
(see \cite{ap:73}), so does $Z$.
The set $K\subset Z$ has property $C$ (as a compact subset of a
finite $C$-space).
Therefore, there exists
a finite $C$-refinement $\{\overline{\gamma _n}\}_{n=1}^{k}$ of
$\{\overline{\omega _n}\cap K\}$ consisting of finite, disjoint and
functionally open in $K$ families. Since $K$ is compact,
we can find finite, disjoint and functionally open in $Z$ families
$\beta _n$, $n=1,..,k$, such that every $\beta _n$ refines
$\overline{\omega _n}$ and $V=\cup\{G\in\beta _n:n=1,..,k\}$ covers
$K$. Next, take an open set $W\subset Z$ such that
$Z\backslash\overline{W}$ is functionally open in $Z$ and
$K\subset W\subset\overline{W}\subset V$. Then $Z\backslash W$, being
closed in $Z$ and disjoint from $K$, is finite-dimensional.
So is $Z\backslash\overline{W}$ as an $F_{\sigma}$-subset of
$Z\backslash W$. Let
$dim Z\backslash\overline{W}=m$. According to \cite[Proposition 2.12]{ag:78},
there exists a finite $C$-refinement $\{\beta _n\}_{n=k+1}^{k+m+2}$
of $\{\overline{\omega _n}\cap (Z\backslash\overline{W})\}_{n=k+1}^{k+m+2}$.
Since $Z\backslash\overline{W}$ is paracompact, we can assume that
$\{\beta _n:n=k+1,..,k+m+2\}$ is a locally finite and functionally open cover
of $Z\backslash\overline{W}$. We finally obtain that
$\{\beta _n\}_{n=1}^{k+m+2}$ is a finite $C$-refinement for the
sequence $\{\overline{\omega _n}\}$ such that $\{\beta _n:n=1,..,k+m+2\}$
is a locally finite and functionally open cover of $Z$. Then
$\{\beta _n\cap X\}_{n=1}^{k+m+1}$ is the required finite $C$-refinement
of $\{\omega _n\}$.

$(b)\Rightarrow (a)$. Let $\omega=\{\omega _n\}$ be a sequence of finite
functionally open covers of $X$. According to our assumption, there is
a finite $C$-refinement $\{\beta _n\}_{n=1}^{k}$ for $\omega$ such
that $\{\beta _n:n=1,..,k\}$ is a locally finite and functionally open
cover of $X$. For fixed $m$ let $\omega _m=\{U_i:i=1,..,s\}$ and
$\beta _m=\{V_{\alpha}:\alpha\in A\}$. We define disjoint sets $A_i\subset A$,
$i=1,2,..,s$, 
by $A_1=\{\alpha:V_{\alpha}\subset U_1\}$ and
$A_{i+1}=\{\alpha\in A\backslash\cup _{j=1}^{i}A_j:V_{\alpha}\subset U_{i+1}\}$.
Since the union of any locally finite family consisting of functionally
open sets is again functionally open, every
$W_i=\cup\{V_{\alpha}:\alpha\in A_i\}$ is functionally open. Moreover,
$\gamma _m=\{W_i:i=1,..,s\}$ is a disjoint family refining $\omega _m$.
Then $\{\gamma _n\}_{n=1}^{k}$ is a finite $C$-refinement of $\omega$
and consists of finite functionally open families.
\end{proof}

\begin{corollary}
Every finite $C$-space has property $C$.
\end{corollary}

Next corollary was actually proved in Theorem 2.4.

\begin{corollary}
If $X$ is paracompact, then $X$ is a finite $C$-space if and only
if $X$ satisfies condition $(K)$ such that $K\subset X$ has
property $C$.
\end{corollary} 

\section*{3. Proof of Theorem 1.1}
\setcounter{section}{3}
\setcounter{theorem}{0}

\bigskip
$(a)\Rightarrow (b)$. Our proof is a slight modification of the proof of \cite[Theorem 2.1]{vu:98}. 
Suppose $X$ is a finite $C$-space and
let $\{\phi _n\}$ be an aspherical sequence of
strongly lsc maps from $X$ into $Y$.
We shall construct by induction the following objects (for 
simplicity, every abstract simplicial complex is identified
with its polyhedron equipped with the Whitehead topology):

\begin{enumerate}
\item[(1)] a sequence of $\Lambda _n$ of pairwise disjoint sets;

\item[(2)] two sequences $\{\omega _n\}$ and $\{\gamma _n\}$
of open locally finite covers of $X$ such that
$\omega _n=\{U_\alpha :\alpha\in\Lambda _n\}$, 
$\gamma _n=\{V_\alpha :\alpha\in\Lambda _n\}$ and
$\overline{U_{\alpha}}\subset V_{\alpha}$ for every 
$\alpha\in\Lambda _n$;

\item[(3)] an increasing sequence of simplicial complexes
$\{\mathcal{K}_n\}$ such that $dim\mathcal{K}_n\leq n-1$ and
$\bigcup _{k=1}^{k=n}\Lambda _k$ is the set of vertices of $\mathcal{K}_n$;

\item[(4)] finite subcomplexes $P_\alpha$ of $\mathcal{K}_n$
for every $\alpha\in\Lambda _n$ and $n\in\N$;

\item[(5)] continuous maps $g_n: \mathcal{K}_n\rightarrow Y$
such that  $g_n$ extends $g_m$ whenever $m<n$ 
and
$g_n(P_\alpha )\subset\phi _n(x)$ for every 
$x\in V_\alpha\in\gamma _n$; 

\item[(6)] if $\bigcap _{j=1}^{j=k}U_{\alpha (j)}\neq\emptyset$,
where $\alpha (j)\in\Lambda _{n(j)}$, $j=1,..,k$, and
$n(1)<....<n(k)$, then the set $\{\alpha (1),..,\alpha (k)\}$ is a
simplex of the complex $P_{\alpha (k)}$.
\end{enumerate}

\medskip
Let describe the first step of the induction. For every 
$x\in X$ fix a point $c(x)\in\phi _1(x)$. Since $\phi _1$
is strongly lsc, we can find neighborhoods $O(x)$ in $X$
such that $z\in O(x)$ yields $c(x)\subset\phi _1(z)$.
Let $\gamma _1=\{V_\alpha :\alpha\in\Lambda _1\}$ 
be a locally finite open refinement of  
of the cover $\{O(x):x\in X\}$ and 
$\omega _1=\{U_\alpha :\alpha\in\Lambda _1\}$ an index
closure refinement of $\gamma _1$. Let $\mathcal{K}_1$
be the $0$-dimensional complex with the set of vertices
$\Lambda _1$ and every $P_\alpha$, $\alpha\in\Lambda _1$
be the subcomplex of $\mathcal{K}_1$ with one vertex 
$\alpha$. Choose the points $x_\alpha\in X$ such that 
$V_\alpha\subset O(x_\alpha )$ and define 
$g_1:\mathcal{K}_1\rightarrow Y$ by 
$g_1(\alpha)=c(x_\alpha)$. 
Following the Uspenskij argumnts from the proof of
\cite[Theorem 2.1]{vu:98} and using that each $\phi _n(x)$ is
$C^{n-1}$-embedded in $\phi _{n+1}(x)$, one can complete the
construction.

\smallskip
As in the proof of \cite[Theorem 2.1]{vu:98},
we can see that property (6) holds.  We set
$\mathcal{K}=\bigcup _{n=1}^{\infty}\mathcal{K}_n$ and
$g:\mathcal{K}\rightarrow Y$ is the continuous map extending
all maps $g_n$.
Since $X$ is a paracompact finite $C$-space, by Theorem 2.4,
there exists a finite sequence
$\{\lambda _k\}_{k=1}^{k=m}$ of open disjoint families such that
$\lambda _k$ refines $\omega _k$ and $\lambda =\bigcup _{k=1}^{k=m}$
is a locally finite cover of $X$. We can suppose that each
$\lambda _k=\{W_\alpha :\alpha\in\Lambda _k\}$ is an index refinement
of $\omega _k$, so
$\lambda =\{W_\alpha :\alpha\in\Gamma\}$, where
$\Gamma =\bigcup _{k=1}^{k=m}\Lambda _k$. Let
$\{h_\alpha :\alpha\in\Gamma\}$
be a locally finite partition of unity subordinated to $\lambda$.
For any $x\in X$ let $s(x)=\{\alpha\in\Gamma :h_\alpha (x)>0\}$.
The sets $s(x)$ are finite and since each $\lambda _k$ is disjoint,
we have $|s(x)\cap\Lambda _k|\leq 1$ for $k\leq m$. Hence
$\displaystyle s(x)=\{\alpha _1,..,\alpha _j\}$ 
such that $\alpha _k\in\Lambda _{i(k)}$ for each $k=1,..,j$
and
$i(1)<..,<i(j)\leq m$.
Observe that $\bigcap _{k=1}^{k=j}\displaystyle U_{\alpha _k}\neq\emptyset$
(it contains $x$). Therefore, by (6), $s(x)$ is a simplex of
$P_{\alpha _{j}}$. Then the formula
$h(x)=\sum\{h_\alpha (x)\alpha :\alpha\in\Gamma\}$ defines a
continuous map $h:X\rightarrow\mathcal{K}_m$ such that
$h(x)\in P_\delta$ for some $\delta\in s(x)$. But $\delta\in s(x)$ implies
$W_\delta\subset U_\delta\subset V_\delta$. Consequently,
by (5), $g_m(h(x))\in g_m(P_\delta )\subset\phi _m(x)$. We finally obtain
that the composition $g_m\circ h$ is a selection for $\phi _m$.

\medskip
$(b)\Rightarrow (c)$ This implication is trivial because every
set-valued map with an open graph is strongly lsc.

\smallskip
$(c)\Rightarrow (a)$ We shall prove that Theorem 1.1(c) implies
Theorem 1.2(b), and then the proof of the present implication will
follow from Theorem 1.2.

We agree the following notations: for any metric space $Y$,
a point $y\in Y$ and a positive number $\delta$ 
the open (resp., closed) ball in $Y$ with center $y$ and radius
$\delta$ is denoted by $B(y,\delta)$ (resp., $\overline{B}(y,\delta)$).
If $Y$ is a normed space, then $d$ stands for the metric 
generated by the norm of $Y$.

Now we need some notations from \cite{gv:99}. Let $Y$ be a metric
space. For every pair of
set-valued maps $\phi\colon X\to 2^Y$ and $\psi\colon X\to 2^Y$
a set-valued map 
$\bigtriangleup _{(\phi ,\psi )}\colon X\times Y\to 2^{\R}\cup\{\emptyset\}$ is
associated such that, for $(x,y)\in X\times Y$, the value
$\bigtriangleup _{(\phi ,\psi )}(x,y)$ consists of all $\delta >0$ for
which there exists a neighborhood
$U_{\delta}$ of $x$ with $B(y,\delta)\cap\phi (z)\neq\emptyset$
and $B(y,\delta)\cap\psi (z)=\emptyset$ whenever $z\in U_{\delta}$. We define
a set-valued map
$\Phi _{(\phi,\psi)}\colon X\to 2^Y\cup\{\emptyset\}$ by 
$\Phi _{(\phi,\psi)}(x)=\{y\in Y: \bigtriangleup _{(\phi ,\psi)} (x,y)\neq\emptyset\}$.
Consider the functions
$u_{(\phi,\psi)}, l_{(\phi,\psi)}\colon G(\Phi _{(\phi,\psi)})\to\R$,
$u_{(\phi,\psi)}(x,y)=\inf\bigtriangleup _{(\phi ,\psi)}(x,y)$ and
$l_{(\phi,\psi)}(x,y)=\sup\bigtriangleup _{(\phi ,\psi)}(x,y)$,
where
$G(\Phi _{(\phi,\psi)})=\{(x,y)\in X\times Y: y\in\Phi _{(\phi,\psi)}(x)\}$
is the graph of $\Phi _{(\phi,\psi)}$.

\begin{lemma}[\cite{gv:99}] Let $Y$ be a metric space, and let
$\phi\colon X\to 2^Y$ and $\psi\colon X\to 2^Y$ be set-valued
maps such that $\phi$ is lsc, $\psi$ has a closed graph and
$\phi (x)\backslash\psi (x)\neq\emptyset$ for every $x\in X$.
Then,

\begin{alphanum}
\item $\phi (x)\backslash\psi (x)\subset\Phi _{(\phi,\psi)}(x)$
for every $x\in X$;

\item The graph of $\Phi _{(\phi,\psi)}$ is open in $X\times Y$;

\item $u_{(\phi,\psi)}$ and $l_{(\phi,\psi)}$ are, respectively,
usc and lsc functions.
\end{alphanum}
\end{lemma}   

\medskip
\begin{proposition}
Theorem $1.1(c)$ yields Theorem $1.2(b)$.
\end{proposition}

\begin{proof}
Suppose that $Y$ is a Banach space,  
$\phi :X\rightarrow {\mathcal F}_c(Y)$ is lsc and
$\{\psi _n\}$ a decreasing sequence
of maps $\psi _n:X\rightarrow {\mathcal F}(Y)$ such that 
each $\psi _n$ has a closed graph and 
$\{\psi _n\}$ is a $Z$-sequence for $\phi$.
Consider the maps $\displaystyle\Phi _n=\Phi_{(\phi,\psi _n)}$.
Observe that each $\phi (x)\backslash\psi _n(x)\ne\emptyset$ because 
$\psi _n(x)\cap\phi (x)$ is a $Z_{n-1}$-set in $\phi (x)$. Then,
by Lemma 3.1, each $\Phi _n$ has an open graph and
$\phi (x)\backslash\psi _n(x)\subset\Phi _n(x)$ for all $x\in X$ and $n$.
Since $\psi _{n+1}(x)\subset\psi _n(x)$,
$\displaystyle\bigtriangleup _{(\phi,\psi _n)}(x,y)$ is a subset of
$\bigtriangleup _{(\phi,\psi _{n+1})}(x,y)$ for every $n$ and
$(x,y)\in X\times Y$. Consquently, $\{\Phi _n\}$ is an increasing
sequence.

\medskip\noindent
{\em Claim. Each $\Phi_n$ is $(n-2)$-aspherical}.

This claim can be proved by making use of the arguments from the
proof of \cite[Claim 1, Proposition 3.1]{gv:99}.

\medskip
We now proceed to the rest of the proof. We already proved that 
$\{\Phi_n\}$ satisfies the hypotheses of Theorem 1.1(c). Hence,
there exists $m\in\N$ such that $\Phi _m$ admits a continuous
selection $g\colon X\to Y$. Relying once again on Lemma 3.1(c)
and the results of \cite{{jd:44},{cd:51},{mk:51}}, there
exists a continuous selection $\delta\colon X\to\R$ for
the map $\displaystyle\bigtriangleup _{(\phi,\psi _m)}(x,g(x))$,
$x\in X$.Then 
$d(g(x),\phi (x))<\delta (x)<d(g(x),\psi _m(x))$ for every $x\in X$.
Therefore, $F(x)=\overline{B(g(x),\delta (x))\cap\phi (x)}$
defines an lsc map $F\colon X\to {\mathcal F}_c(Y)$. 
Hence, applying 
Michael's convex-valued selection theorem \cite[Theorem 3.2"]{m1:56}, we 
get a continuous selection $f$ for $F$.    
Since $F(x)\subset\phi (x)\backslash\psi _m(x)$,
$f$ is as required.
\end{proof}

\section*{4. Proof of Theorem 1.2}
\setcounter{section}{4}
\setcounter{theorem}{0}

\bigskip
First, note that implication $(a)\Rightarrow (b)$ follows
from Theorem 1.1, implication $(a)\Rightarrow (c)$
and Proposition 3.2.
Implication $(b)\Rightarrow (c)$
is trivial, so it remains only to prove the implication $(c)\Rightarrow (a)$.

\medskip\noindent
$(c)\Rightarrow (a)$.
Take a sequence $\{\omega_n:n=1,2,\dots\}$ of finite open covers of $X$.
By definition, we must prove that there exists a finite sequence
$\{\gamma _n:n=1,2,..,m\}$ of disjoint open families in $X$
such that each $\gamma _n$ refines $\omega_n$ and the union
$\bigcup\{\gamma_n:n=1,2,..,m\}$ is a cover of $X$. To this end, we
proceed just like in
\cite[Theorem 1.3, implication $S_2\Rightarrow C$]{vu:98},
(see also \cite[Theorem 1.1, implication
$(c)\Rightarrow (a)$]{gv:99}) with a few modifications.
Note that, for every $n$, there exists an $\omega _n$-map
$f_n\colon X\to Q_n$ of $X$ to a compact polyhedron $Q_n$. Considering
$Q_n$ as a subspace of a Euclidean space ${\R}^{k(n)}$  
and multiplying,
if necessary, $f_n$ by a constant, we may assume that $Q_n$ is the cube
$[0,k(n)]^{k(n)}$ and that 

\medskip\noindent
(7) \hbox{}~~~~~~$\{f_n^{-1}(U): U\ \hbox{is an open ball of radius 2 in}\ (Q_n,d_n)\}$
is an open cover of $X$ refining $\omega_n$.

\medskip\noindent
Here we consider the metric $d_n((x_j),(y_j))=\max |x_j-y_j|$ on $\R^{k(n)}$.
Let $A_n$ be the union of all sets $B(j,s)=\{(x_m)\in Q_n:x_j=s\}$, where
$j,s$ are natural numbers with $1\leq j\leq k(n)$ and $1\leq s\leq k(n)-1$.
Obviously, $Q=\prod _{n=1}^{\infty}Q_n$ is homeomorphic
to the Hilbert cube. Next, let $Z$ be the linear subset of
$E=\prod _{n=1}^{\infty}\R^{k(n)}$ defined as

\[
Z=\left\{(y_{k(n)})\in E:\sum_{n=1}^\infty
\frac{\|y_{k(n)}\|}{2^{k(n)}}<+\infty\right\}.
\]
The space $Z$ is equipped with the norm 
\[
\|(y_{k(n)})\|_\omega=\sum_{n=1}^\infty \frac{\|y_{k(n)}\|}{2^{k(n)}},
\]
where $\|y_{k(n)}\|=\max |y_{k(n),j}|$, $y_{k(n),j}$
being the $j$th coordinate of $y_{n(k)}$,
is the norm on $\R^{k(n)}$.
Observe that the
topology of $(Z,\|.\|_\omega)$ coincides with the topology of $Z$
as a subspace of $E$, and let $Y$ be the completion of $Z$.

It suffices to find a finite sequence
$g_n\colon X\to Q_n$, $n=1,..,m$, of continuous maps such that

\medskip\noindent
(8) \hbox{}~~~$d_n(f_n(x),g_n(x))\leq 1$ for all $x\in X$ and $n\leq m$;

\smallskip\noindent
(9) \hbox{}~~~For every $x\in X$ there is $n\leq m$ with $g_n(x)\not\in A_n$.

\medskip
Indeed, let $\lambda _n$ be the family of all components of
$Q_n\backslash A_n$ which is disjoint and consists of open sets in $Q_n$
with diameter $\leq 1$. Then (8) yields that each of the 
disjoint families
$\gamma _n=g_n^{-1}(\lambda _n)$ refines $\omega _n$ 
(see \cite[proof of Theorem 1.3]{vu:98}), and (9)
implies that $\bigcup _{n=1}^{n=m}\gamma _n$ covers
$X$. Let us observe also that condition (9) is equivalent
to the existence of a continuous map
$G:X\rightarrow Q(m)=\prod _{n=1}^{m}Q_n$ avoiding 
the set $A(m)=\prod _{n=1}^{m}A_n$. 

We shall prove that, under the hypotheses of (c), such
$g_n$'s exist. For every $n$ consider the map 
$\Phi_n:X\to {\mathcal F}_c(\R^{k(n)})$ defined by
\[
\Phi_n(x)=\{y\in Q_n:d_n(f_n(x),y)\leq 1\},\quad x\in X 
\]
and define $\Phi:X\to{\mathcal F}_c(Y)$ by
\[
\Phi(x)=\prod\{\Phi_n(x):n=1,2,\dots\},\quad x\in X.
\]
It is easily seen that $\Phi$ is lsc and   
$\Phi(x)\subset Q\subset Z$ for every $x\in X$.

\medskip
We proceed to the final
step in this proof. Thus, we have a Banach space $Y$ and an lsc
map $\Phi:X\to {\mathcal F}_c(Y)$. 
Further, if $A=\prod _{n=1}^{\infty}A_n\subset Q$, then for every $x\in X$,
$\Phi (x)\cap A$ is the product $\prod _{n=1}^{\infty}A_n\cap\Phi _n(x)$.
Since each $A_n\cap\Phi _n(x)$ is closed and nowhere dense in $\Phi _n(x)$,
by \cite[Corollary 2]{bt:99},
$\prod _{n=1}^{n=k}A_n\cap\Phi _n(x)$ is a
$Z_{k-1}$-set in $\prod _{n=1}^{n=k}\Phi _n(x)$ for
every $k$. This yields that $\Phi (x)\cap F_k$ is a $Z_{k-1}$-set in $\Phi (x)$,
$x\in X$, where $F_k=A(k)\times\prod _{n=k+1}^{\infty}Q_n$. 
Therefore the map $\Phi$ and the decreasing sequence $\{F_k\}$
satisfy the hypotheses of (c). Hence, there is $m$ and a 
continuous selection $G_1$ for $\Phi$ avoiding $F_m$. 
Let $\pi _m$ be the projection from 
$Q$ onto $Q(m)$ and $G=\pi _m\circ G_1$. 
Then $G=(g_1,..,g_m)$, where each $g_n$ is
a continuous map from $X$ into $Q_n$. It is easily seen that the maps $g_n$
satisfy (8) and (9). Hence, the proof is complete.

\section*{5. Some corollaries and applications}
\setcounter{section}{5}
\setcounter{theorem}{0}

Let $(Y_n,\|.\|_n)$ be a sequence of normed spaces. We denote
by $Y$ the product $\prod _{n=1}^{\infty}Y_n$ and for
every $m$ let $\pi _m\colon Y\to Y(m)$ be the natural projection, 
where $Y(m)=\prod _{n=1}^{m}Y_n$.  
We say that a sequence of maps $f_n\colon X\to Y_n$ is 
finitely removable from a set $A\subset Y$ if for every
sequence $\{\epsilon _n\}$ of positive real numbers there
exists a finite sequence of maps $g_n\colon X\to Y_n$,  
$n=1,..,m$, such that $g(X)\cap\pi _m(A)=\emptyset$
and
$\|(f_n(x)-g_n(x)\|_n\leq\epsilon _n$ for every
$n$ and $x\in X$,
where $g=(g_n)\colon X\to Y(m)$.
A sequence of maps $f_n\colon X\to Y_n$
is called uniformly bounded if there exists a common
bound for all diameters $diam(f_n(X))$, $n\in\N$.

Our first application of Theorem 1.2 is the following characterization
of finite $C$-spaces:

\begin{proposition}
A space $X$ is a finite $C$-space if and only if it satisfies
the following condition: 
for any sequence of Banach spaces $Y_n$ and closed nowhere
dense subsets $A_n\subset Y_n$,
every uniformly bounded sequence of maps $f_n\colon X\to Y_n$
is finitely removable from the set $A=\prod _{n=1}^{\infty}A_n$.
\end{proposition}

\begin{proof} Suppose $X$ satisfies the condition from
Proposition 5.1.
To show that $X$ is a finite $C$-space, take a sequence
$\{\omega _n\}$ of finite functionally open covers of $X$. 
Proceeding just like in the proof of Theorem 1.2, implication
$(c)\Rightarrow (a)$, we obtain $\omega _n$-maps 
$f_n\colon X\to Q_n\subset\R^{k(n)}$, where now $Q_n=I^{k(n}$,  
$I=[1,2]$. Let $K_n=[0,3]^{n(k)}$.
Since $f_n\colon X\to K_n$ is an $\omega _n$-map, there exists
positive $\epsilon _n\leq 1$ such that if $U$ is an open ball in 
$K_n$ with radius $<2\epsilon _n$, then $f_n^{-1}(U)$
is contained in an element of $\omega _n$.
Next, take hyperplanes
$B(j,s)=\{(x_i)\in\R^{k(n)}:x_j=s\}$, $1\leq j\leq k(n)$ and
$0<s<3$, such that their union $A_n$ partitions $K_n$ into
open (in $K_n$) cubes, each having a diameter $<\epsilon _n$.
Since $\{f_n\}$ is uniformly bounded and every
$A_n\subset\R^{k(n)}$ is closed nowhere dense, 
$\{f_n\}$ is finitely removable
from the set $A=\prod _{n=1}^{\infty}A_n$. So, there exists a
finite sequence of maps $g_n\colon X\to\R^{k(n)}$, $n=1,..,m$
such that each $g_n$ is $\epsilon _n$-close to $f_n$ and
the map $g=(g_n)\colon X\to\prod _{n=1}^{m}\R^{k(n)}$ avoids
the set $A(m)=\prod _{n=1}^{m}A_n$. Observe that all $g_n$ map $X$
into $K_n$ (because $g_n$ is $\epsilon _n$-close to $f_n$) and
the family $\lambda _n$ of all components of $K_n\backslash A_n$
is disjoint and consists of open sets in $K_n$ with diameter
$<\epsilon _n$. Consequently, $\gamma _n=g_n^{-1}(\lambda _n)$
is a disjoint and functionally open in $X$ family which refine
$\omega _n$. It remains only to note that, since the map $g$
avoids the set $A(m)$, $\{\gamma _n:n=1,..,m\}$ covers $X$.
Hence, $\{\gamma _n\}$ is a finite $C$-refinement for $\{\omega _n\}$. 

To prove the inverse implication, first let show
that $X$ can be supposed to be paracompact. Indeed, let
the maps $f_n\colon X\to Y_n$ and the sets $A_n\subset Y_n$
be as in Proposition 5.1. Consider the \v Cech-Stone extension
$\beta f\colon\beta X\to\beta Y$ of $f=(f_n)$ and let
$H=(\beta f)^{-1}(Y)$. Obviously, $H$ is paracompact (as
a perfect preimage of $Y$).
Denote by $h=(h_n)\colon H\to Y$ the
restriction $(\beta f)|H$. Then $h_n\colon H\to Y_n$ is
uniformply bounded and, by Proposition 2.2, $H$ is a finite
$C$-space. Hence, if $\{h_n\}$ is finitely removable from $A$,
so is $\{f_n\}$.
Therefore, we can assume that $X$ is paracompact.

We consider the normed space $Z$ defined by

\[
Z=\left\{(y_n)\in Y:\sum_{n=1}^\infty
\frac{\|y_n\|_n}{2^n}<+\infty\right\},
\]
where
\[
\|(y_n)\|=\sum_{n=1}^\infty \frac{\|y_n\|_n}{2^n}
\]
is the norm in $Z$. 
Let $E$ be the completion of $Z$ and,
for a given bounded sequence $\{\epsilon _n\}$ of positive numbers,
define the lsc map $\phi\colon X\to {\mathcal F}_c(Z)$,
$\phi (x)=\{(y_n)\in Y:\|y_n-f_n(x)\|_n\leq\epsilon _n, n\in\N\}$.
This definition is correct because $\{f_n\}$ is uniformly bounded.
The norm topology of $Z$ coincides with the topology inherited
from $Y$ and the values of $\phi$ are norm-complete in $Z$
(see
\cite[proof of Theorem 1.1, implication $(c)\Rightarrow (a)$]{gv:99}
for a similar proof).
So,
$\phi$ is a map from $X$ into ${\mathcal F}_c(E)$.
Let
$D_n=A^n\times\prod _{k=n+1}^{\infty}Y_k$ and $F_n$ be the closure
in $E$ of $D_n\cap Z$, $n\in\N$.
Then, using a result of \cite{bt:99}, we can show that
each $F_n\cap\phi (x)$ is a $Z_{n-1}$-set in $\phi (x)$. 
Since $\{F_n\}$
is a decreasing sequence, by Theorem 1.2(c), there exists
$m\in\N$ and a continuous map $G\colon X\to E$ such that
$G(x)\in\phi (x)\backslash F_m$ for every $x\in X$. Finally,
the map $g=\pi _m\circ G$ is generated by a sequence
$g_n\colon X\to Y_n$, $n=1,..,m$ which
witnesses that $\{f_n\}$ is finitely removable from the set $A$.
\end{proof}

Next corollary, which follows from the proof of Proposition 5.1,
is useful for applications.

\begin{corollary} 
$X$ is a finite $C$-space iff for any sequence $\{Q_n\}$ of
finite-dimensional cubes
and closed nowhere dense sets $A_n\subset Q_n$, every sequence
of maps
$f_n\colon X\to Q_n$ is finitely removable
from the set $A=\prod _{n=1}^{\infty}A_n$.
\end{corollary}

Following F. Ancel \cite{fa:85}, we say that a map $f\colon X\to Y$ is
approximately invertible if there exists a $C^*$-embedding
$i\colon X\to Z$ into a space $Z$ such that for every collection
${\mathcal W}$ of open subsets of $Z$ which is refined by
$\{i(f^{-1}(y)):y\in Y\}$, there is a map $g\colon Y\to Z$
with $g\circ f$ is ${\mathcal W}$-close to $i$ (i.e. the family
$\{g(f(x)),i(x):x\in X\}$ refines $\mathcal W$). F. Ancel
\cite{fa:85} proved that every approximately invertible 
surjection between metric spaces having compact fibers preserves
property $C$. 
The analogue of Ancel's result for finite $C$-spaces also holds.

\begin{corollary}
Any approximately invertible surjection 
preserves finite $C$-space property.
\end{corollary}

\begin{proof} Suppose $p\colon X\to Y$ is approximately invertible
and surjective, where $X$ 
is a finite $C$-space. Let $\{Q_n\}$ be a sequence of 
finite-dimensional cubes, $A_n\subset Q_n$ closed nowhere dense
subsets and
$f_n\colon Y\to Q_n$ arbitrary sequence of maps.
By Corollary 5.2,
if suffices to show that $\{f_n\}$ is finitely removable from the set
$A=\prod _{n=1}^{\infty}A_n$.
To this end, let $d_n$ be a metric on $Q_n$ and 
$\{\epsilon _n\}$ a sequence of positive numbers.
Since $X$ is a finite $C$-space, according to
Corollary 5.2,
$\{f_n\circ p\}$ is finitely removable from
$A$. So, there exists a number $m\in\N$ and a map
$g=(g_n)\colon X\to Q(m)=\prod _{n=1}^{m}Q_n$ such that
$g(X)\cap A(m) =\emptyset$ and
$d_n(f_n(p(x)),g_n(x))\leq 4^{-1}\epsilon _n$ for every $x\in X$ and
$n\leq m$, where $A(m)=\prod _{n=1}^{m}A_n$.
Now, let $i\colon X\to Z$ be a $C^*$-embedding into a space $Z$ which
witnesses the approximate invertibility of the map $p$.
We identify $i$ with the identity $id_X$ and
take an extension $\overline{g}\colon Z\to Q(m)$ of $g$. Let
$U=\{z\in Z:\overline{g}(z)\not\in A(m)\}$. For every $y\in Y$ 
let
$W(y)=\{z\in U:d_n(\overline{g}_n(z),f_n(y))<3^{-1}\epsilon _n,n\leq m\}$
and ${\mathcal W}=\{W(y):y\in Y\}$.
Obviously, each $W_y$ is open in $Z$ and contains $p^{-1}(y)$.
So, there is a map
$q\colon Y\to Z$ such that $q\circ p$ and $id_X$ are
${\mathcal W}$-close. This implies that $q(Y)\subset U$.  
Define now $h_n\colon Y\to Q_n$ by 
$h_n=\overline{g}_n\circ q$, $n\leq m$. Then
$h=(h_n)\colon Y\to Q(m)$ avoids the set $A(m)$.

It remains to show that each $h_n$ is $\epsilon _n$-
close to $f_n$, $n\leq m$.
For any $y\in Y$ choose 
$x\in p^{-1}(y)$. Then $x$ and $q(p(x))=q(y)$ belong
to $W(y')$ for some $y'\in Y$ because $q\circ p$ and $id_X$
are $\mathcal W$-close. Hence, 
$d_n(\overline{g}_n(x),\overline{g}_n(q(y)))<2.3^{-1}\epsilon _n$
for every $n\leq m$. On the other hand, since 
$\overline{g}_n(x)=g(x)$, we have
$d_n(f_n(y),\overline{g}_n(x))\leq 4^{-1}\epsilon _n$.
Finally, $d_n(f_n(y),\overline{g}_n(q(y)))\leq\epsilon _n$.
\end{proof}

Y. Hattori and K. Yamada \cite{hy:89} proved that if 
$f\colon X\to Y$ is a closed map from a paracompact space $X$ onto
a $C$-space $Y$ and all fibers $f^{-1}(y)$ have property $C$, then
$X$ is a $C$-space. Similar result for finite 
$C$-spaces is not valid. Indeed, let $X$ be a disjoint union
of $n$-dimensional cubes $Q_n$ and for every $n\in\N$ denote by
$f_n$ the projection of $Q_n$ onto $I_n=[0,1]$. Let $Y$ be the
disjoint union of all $I_n$ and define $f\colon X\to Y$ as
$f|Q_n=f_n$. Then, $f$ is a perfect map with finite-dimensional
fibers and $Y$ is a finite $C$-space ($dim Y=1$). Since every
paracompact finite $C$-space satisfies condition
$(K)$ (see Corollary 2.6), $X$ is not a finite $C$-space.
But still there is an analogue of the Hattori and Yamada
result.

\begin{proposition}
Let $f\colon X\to Y$ be a closed surjection such that $X$ is normal, 
$Y$ paracompact and $dim f^{-1}(y)\leq k$ for every $y\in Y$.
If $Y$ is a finite $C$-space, so is $X$.
\end{proposition}

\begin{proof}
Let $T=(\beta f)^{-1}(Y)$. Since $f$ is closed and $X$ normal, the closure in $\beta X$
of each fiber $f^{-1}(y)$, $y\in Y$, is $(\beta f)^{-1}(y)$. Hence,
$dim (\beta f)^{-1}(y)\leq k$ for any $y\in Y$. Because $X$ and $T$ have
the same \v Cech-Stone compactification, they are simultaneously finite
$C$-spaces.
Therefore, we can suppose that $f$ is perfect and $X$ is paracompact
(as a perfect preimage of $Y$).

Since $Y$ is a paracompact finite $C$-space, 
it satisfies condition $(K)$ such that $K\subset Y$ has
property $C$. Then, by the mentioned result of Hattori and
Yamada, $H=f^{-1}(K)$ is a compact $C$-space. Hence, according
to Corollary 2.6, it suffices to show that every closed
$F\subset X$, which is disjoint from $H$, is finite-dimensional. 
To this end, observe that $f(F)\subset Y$ is closed and
disjoint from $K$, so $f(F)$ is finite-dimensional. Then, by the
Hurewicz formula (applied to the map $f|F\colon F\to f(F)$, 
see \cite{es:62}), $dimF\leq dimf(F)+k$.  
\end{proof}

\section*{6. Factorization theorems}
\setcounter{section}{6}
\setcounter{theorem}{0}

Let $f\colon X\to Y$ be a continuous surjection. We say that $(Z,h,g)$
is a factorization for $f$ if $h$ is a continuous map from $X$
to the space $Z$ and $g$ is a continuous map from $Z$ onto $Y$
such that $f=g\circ h$ and $w(Z)\leq w(Y)$. Here $w(Y)$ is the
topological weight of $Y$.
V. Chatyrko \cite{cha:99} proved that if $f\colon X\to Y$ is a map
between compact spaces and $X$ has property $C$, then there is a
factorization $(Z,h,g)$ such that $Z$ is a compact $C$-space and
$dim_CZ\leq dim_CX$, where $dim_C$ is the Borst transfinite
dimension.
Another factorization theorem for $C$-spaces was proved
by M. Levin, L. Rubin and P. Schapiro \cite{lrs:98} which 
implies that if $f\colon X\to Y$ is a map between compact spaces
and $\{D_n\}$ is a sequence of compact subsets of $X$ with property
$C$, then $f$ admtis a factorization $(Z,h,g)$ such that
each $h(D_n)$ has property $C$.
In this section
we are presenting some factorizations for $C$-spaces
in the special case when $Y$ is a metric space.

\begin{proposition}
Let $\{K_n\}$ be a sequence of compact $C$-subsets of a space $X$. 
Then every map $f\colon X\to Y$ into a metrizable space $Y$
admits a factorization $(Z,h,g)$
such that $Z$ is metrizable and all sets $h(K_n)$ have property $C$. 
\end{proposition}

\begin{proof}
Recall that a map $p\colon Y\to S$ between metrizable spaces is
strongly $0$-dimensional if there is a metric on $Y$ generating
its topology such that for every $\epsilon >0$ and every
$z\in p(Y)$ there exists an open neighborhood $U$ of $z$ with
$p^{-1}(U)$ being the union of disjoint open sets of diameter
$<\epsilon$. It is well known that every metrizable space
admits a strongly $0$-dimensional map into Hilbert cube $Q$.
So, we can take a strongly $0$-dimensional map $p\colon Y\to Q$
and let $\overline{f}\colon\beta X\to Q$ be the \v Cech-Stone
extension of $p\circ f$. By the Levin, Rubin and Schapiro result
mentioned above, there is a factorization $(Z_0,h_0,g_0)$
for $\overline{f}$
such that all sets $h_0(K_n)$ have property $C$. Let
$g\colon Z\to Y$ and $q\colon Z\to h_0(X)$ be the pullback of
$p$ and $g_0|h_0(X)$ respectively, and $h\colon X\to Z$,
$h(x)=(h_0(x),f(x))$. Then $Z$ is a metrizable space of the same
weight as $Y$ and $Z$ admits a compatible metric such that $q$
is strongly $0$-dimensional. Since all fibers of $q$
are $0$-dimensional, each
restriction $q_n=q|h(K_n)$ is a map from $h(K_n)$ onto
$h_0(K_n)$ with $0$-dimensional fibers. Hence, by \cite{hy:89},
$h(K_n)$ has property $C$.
\end{proof}

\begin{proposition}
Every map $f\colon X\to Y$ from a finite $C$-space $X$ onto
metrizable $Y$ admits a factorization $(Z,h,g)$ such that
$Z$ is a metrizable finite $C$-space.
\end{proposition}

\begin{proof}
Let $P=(\beta f)^{-1}(Y)$ and $\overline{f}$ be the restriction $(\beta f)|P$.
Then $P$ is a paracompact finite
$C$-space, so it satisfies condition $(K)$ such that 
$K\subset P$ is a compact $C$-space.
Let $K_Y=\overline{f}(K)$ and $K_X=(\overline{f})^{-1}(K_Y)$.
Obviously, $K_X\subset P$ is compact and contains $K$.
Next, if $d$ is a metric on $Y$,
let $H_n=\{y\in Y:d(y,K_Y)\geq n^{-1}\}$ and
$F_n=(\overline{f})^{-1}(H_n)$.
Since each $F_n$ is closed in $P$ and disjoint from $K_X$, it
is finite-dimensional.
Because $\overline{f}$ is a closed map and
any closed set in $Y$ disjoint from $K_Y$ is contained in some
$H_n$,  we have that

\medskip\noindent
(10) every closed set in $P$ disjoint from $K_X$ is contained in some
$F_n$.

\medskip\noindent
Observe that $K_X$ has property $C$ as a compact set in $P$.
Hence, by Proposition 6.1, $\overline{f}$ admits a factorization
$(M,h_1,g_1)$ such that $M$ is metrizable and the set
$K_M=h_1(K_X)$ has property $C$. 
Now, we can apply Pasynkov's factorization theorem \cite{bp:68} to
obtain a factorization $(G,h_2,g_2)$ for the map $h_1$ such that
$G$ is metrizable and $dim\overline{h_2(F_n)}\leq dim F_n$ for every
$n$. But each $h_2(F_n)$ is closed in $G$, so $G$ is the union of
the closed sets $K_G=h_2(K_X)$ and $G_n=h_2(F_n)$ with each $G_n$
being finite-dimensional and disjoint from $K_G$. It follows
from (10) that 
every closed subset of $G$ disjoint from $K_G$ is contained in some
$G_n$. Finally, attaching to $G$ the space $K_M$ by the map
$g_2|K_G$, we obtain the metrizable space $Z$. Obviously, $K_M$
is homeomorhic to a subset $K_Z\subset Z$ and $Z\backslash K_Z$
is the set $G\backslash K_G$. Therefore, $Z$, being the union of
all sets $G_n$, $n\in\N$, and $K_Z$, satisfies condition $(K)$
with the compact $C$-set $K_Z$. Hence, $Z$ is a finite $C$-space.
It remains only to define the maps
$h\colon X\to Z$ and $g\colon Z\to Y$ such that $f=g\circ h$.
Let $q\colon G\to Z$ be the natural quotient map. Identifying
$K_Z$ with $K_M$ we set $h(x)=q(h_2(x))$, $x\in X$, and 
$g(z)=g_1(g_2(z))$ provided $z\in Z\backslash K_Z$ and $g(z)=g_1(z)$
if $z\in K_Z$.
\end{proof}

\begin{corollary}
Every metrizable finite $C$-space has a completion which is also
a finite $C$-space.
\end{corollary}

\begin{proof}
Suppose $X$ is a metrizable finite $C$-space and let $Y$ be any
completion of $X$. Extend the inclusion $X\subset Y$ to a map
$p\colon\beta X\to\beta Y$ and denote $G=p^{-1}(Y)$ and 
$f=p|G\colon G\to Y$. Then $G$ is a finite $C$-space, so, by
Proposition 6.2, $f$ admits a factorization $(Z,h,g)$ with $Z$
a metrizable finite $C$-space. Obviously, $h$ embeds $X$ into
$Z$. Finally, since $g$ is a perfect map from
$Z$ onto $Y$ and $Y$ is \v Cech-complete, so is $Z$. 
\end{proof}

\bibliographystyle{amsplain}

\bibliographystyle{amsplain}
\bibliography{triquot}

\end{document}